%% file: EllipticDoubleZetaValues.tex
\documentclass[11pt,a4paper]{article}

\author{Nils Matthes}

\date{}

\makeindex

\include{template}

\RequirePackage{verbatim}

\makeatletter
\newwrite\bibinl@out
\newenvironment{bibtex}[1][\jobname]{%
 \immediate\openout\bibinl@out #1.bib
 \immediate\write\bibinl@out{\@percentchar generated from `\jobname' starting line \the
\inputlineno^^J}%
 \def\verbatim@processline{\immediate\write\bibinl@out{\the\verbatim@line}}%
 \@bsphack\let\do\@makeother\dospecials\catcode`\^^M\active\verbatim@start
}%
{\immediate\closeout\bibinl@out\@esphack}
\makeatother

\begin{document}
\title{Elliptic Double Zeta Values}

\maketitle

\abstract{We study an elliptic analogue of multiple zeta values, the \textit{elliptic multiple zeta values} of Enriquez, which are the coefficients of the \textit{elliptic KZB associator}. Originally defined by iterated integrals on a once-punctured complex elliptic curve, it turns out that they can also be expressed as certain linear combinations of indefinite iterated integrals of Eisenstein series and multiple zeta values.
In this paper, we prove that the $\Q$-span of these elliptic multiple zeta values forms a $\Q$-algebra, which is naturally filtered by the length and is conjecturally graded by the weight. Our main result is a proof of a formula for the number of $\Q$-linearly independent elliptic multiple zeta values of lengths one and two for arbitrary weight.}

\section{Introduction}

1.1. {\bf Multiple zeta values.}
Multiple zeta values are real numbers, defined for a collection $n_1,\dots,n_r$ of positive integers, with $n_r > 1$, by the nested sum
\[
\zeta(n_1,\dots,n_r)=\sum_{0<k_1<k_2<\dots <k_r}\frac{1}{k^{n_1}_1k^{n_2}_2\dots k^{n_r}_r}.
\]
The \textit{weight} is the quantity $n_1+\dots+n_r$ and the \textit{depth} is the number $r$. They were first studied by Euler in the 18th century, and are known to appear in many areas of mathematics and of mathematical physics.

Multiple zeta values satisfy very many $\Q$-linear relations, and it is an important problem to classify those. A conjecture, due to Ihara, Kaneko and Zagier~\cite{IKZ}, states that all such relations are obtained from the \textit{extended double-shuffle relations}, which, roughly speaking, arise from rewriting the product of two multiple zeta values in two different ways as $\Q$-linear combinations of other multiple zeta values. There are also conjectural formulas, due to Broadhurst and Kreimer~\cite{BK}, for the dimensions of the $\Q$-vector space $\gr^{\cD}_d\cZ_N$ of multiple zeta values of a fixed weight and depth. At present however, neither conjecture has been resolved.
\bigskip

\noindent
1.2. {\bf Elliptic multiple zeta values.}
In this paper, we study the \textit{elliptic multiple zeta values} of Enriquez~\cite{E}. These are defined by iterated integrals on a once-punctured complex elliptic curve, and they come in two types, namely \textit{A-elliptic} and \textit{B-elliptic} multiple zeta values, each of which corresponds to one of the two natural homology cycles, $\alpha$ and $\beta$, on an elliptic curve.

Elliptic multiple zeta values were first constructed in the general framework of an elliptic extension of Drinfel'd's version of Grothendieck-Teichm\"uller theory~\cite{E2}, where they arise as coefficients of the \textit{elliptic KZB associator}, which describes the monodromy of the Knizhnik-Zamolodchikov-Bernard (KZB) connection~\cite{CEE,LR}. A second description of elliptic multiple zeta values is available as values of homotopy invariant iterated integrals on a once-punctured elliptic curve~\cite{BL}. The equivalence of these two approaches essentially comes down to the fact, proved in~\cite{BL}, that all such iterated integrals can be constructed from solutions of the KZB equation (cf.~Section~3.2 for more details). This two-fold description of elliptic multiple zeta values is analogous to the fact that multiple zeta values can be defined as coefficients of the \textit{Drinfel'd associator}, which describes the monodromy of the Knizhnik-Zamolodchikov connection~\cite{Dr,LM}, or equivalently as values of homotopy invariant iterated integrals on the projective line minus $0$, $1$ and $\infty$.

Elliptic multiple zeta values are functions on the upper half-plane $\H$. To a point $\tau \in \H$, one associates the framed, once-punctured elliptic curve $E^{\times}_{\tau}=\C/(\Z+\Z\tau) \setminus \{0\}$. In~\cite{BL}, Section~3, it is described how one can associate to $E^{\times}_{\tau}$ an infinite family 
\begin{equation} \label{eqn:diffforms}
\{\nu,\omega^{(k)}\}_{k \geq 0} \subset \cE^1(E^{\times}_{\tau})
\end{equation}
of differential one-forms on $E^{\times}_{\tau}$, which are constructed from a real analytic variant of a classical Kronecker series~\cite{W}, which trivializes the elliptic KZB connection. Elliptic multiple zeta values are now defined as the values of iterated integrals of the forms $\nu,\omega^{(k)}$ along certain paths on the elliptic curve. There are two natural choices of such paths, namely the images $\alpha$ and $\beta$ under the canonical projection $\C \rightarrow E_{\tau}$ of the straight line paths from $0$ to $1$ and from $0$ to $\tau$ in $\C$.\footnote{A technical complication arises from the fact that the start- and endpoints of $\alpha$ and $\beta$ are not contained in the punctured elliptic curve $E^{\times}_{\tau}$. But this can be dealt with using the concept of tangential base point (cf.~\cite{Del}, \S~15).}
The paths $\alpha,\beta$ on $E_{\tau}$ each correspond to different families of elliptic multiple zeta values, A-elliptic and B-elliptic multiple zeta values, and each of them gives rise to a $\Q$-algebra $\eZ^A$ and $\eZ^B$. However, the existence of a modular transformation formula between the two types of elliptic multiple zeta values (cf.~\cite{E}, Section~2.5) allows to transfer a lot of information from A-elliptic to B-elliptic multiple zeta values and vice versa. In particular, there is an explicit formula for B-elliptic multiple zeta values in terms of A-elliptic multiple zeta values (cf.~\cite{E}, eq.~(26)). Another thing is that A-elliptic multiple zeta values are slightly simpler objects than B-elliptic multiple zeta values: the former have a Fourier expansion in $q=e^{2\pi i\tau}$, whereas B-elliptic multiple zeta values only have an expansion in $q$ and $\log(q)$. It is mainly for these two reasons that in this paper we only deal with A-elliptic multiple zeta values.
\bigskip

\noindent
1.3. {\bf Elliptic multiple zeta values and iterated Eisenstein integrals.}
It is convenient to consider the generating series of elliptic multiple zeta values, which is a formal series in two non-commuting variables $x,y$. As mentioned above, this is Enriquez's elliptic KZB associator. The elliptic KZB associator is essentially a tuple $(A(\tau),B(\tau))$ of two formal power series in the variables $x,y$, which depend holomorphically on the natural coordinate $\tau$ on the upper half-plane. One of its most important properties is that both $A$ and $B$ satisfy the following differential equation (cf.~\cite{E2}, Section~5.4), which was already implicit in~\cite{CEE,LR}
\begin{equation} \label{eqn:diffass}
2\pi i\frac{\partial}{ \partial \tau}g(\tau)=\left(-\sum_{ k\geq 0}(2k-1)G_{2k}(\tau)\varepsilon_{2k}\right)g(\tau).
\end{equation}
Here, $G_{2k}$ denotes the usual Eisenstein series (with $G_0 \equiv -1$) and the $\varepsilon_{2k}$ are special derivations of the free Lie algebra on the set $\{x,y\}$~\cite{CEE,Pol}.\footnote{In~\cite{CEE}, these derivations are denoted by $\delta_{2k}$.}
The second important property of the KZB associator is that in the limit $\tau \to i\infty$, it degenerates to a product of Drinfel'd associators (cf.~\cite{E}, eq.~(1.2.4)). Since the coefficients of the Drinfel'd associator are given by multiple zeta values, the upshot of this discussion is that elliptic multiple zeta values are linear combinations, with multiple zeta values as the coefficients, of iterated integrals of Eisenstein series, or \textit{iterated Eisenstein integrals} for short~\cite{Man,MMV}. It is a natural question to determine precisely which linear combinations of multiple zeta values and iterated Eisenstein integrals appear. This question is in fact the main motivation for the present paper.

From \eqref{eqn:diffass}, one sees that the linear combinations of iterated Eisenstein integrals which appear as elliptic multiple zeta values are constrained by relations between the derivations $\varepsilon_{2k}$~\cite{Pol}. More precisely, every linear relation between commutators of the $\varepsilon_{2k}$ gives rise to a linear constraint on the combinations of iterated Eisenstein integrals appearing in elliptic multiple zeta values. We refer to~\cite{BMS} for more details on the relation between the Lie algebra of derivations $\fu^{geom}$ spanned by the $\varepsilon_{2k}$ and elliptic multiple zeta values.

On the other hand, while $\fu^{geom}$ controls the ``iterated Eisenstein integral portion'' of elliptic multiple zeta values, it does not seem to yield information about the ``multiple zeta value portion'' of elliptic multiple zeta values. Moreover, passing to the Lie algebra $\fu^{geom}$ eliminates non-trivial products of elliptic multiple zeta values. Thus, the main difference between~\cite{BMS} and the present paper is that in this paper, information about which multiple zeta values and which products of elliptic multiple zeta values occur is retained.

\bigskip

\noindent
1.4. {\bf The main result}
The principal algebraic object studied in this paper is the $\Q$-algebra of A-elliptic multiple zeta values $\eZ^A$. The notions of weight and depth of multiple zeta values have analogues for elliptic multiple zeta values, namely the \textit{weight} and the \textit{length}. The latter defines an ascending filtration $\cL_{\bullet}\eZ^A_N$ on $\eZ^A_N$, where $\eZ^A_N$ denotes the $\Q$-vector space spanned by elliptic multiple zeta values of weight $N$. Denote by $\gr^{\cL}_{\bullet}$ the associated graded. Since there are only finitely many elliptic multiple zeta values of a fixed weight and length, $\gr^{\cL}_l\eZ^A_N$ is finite-dimensional for every $l$ and $N$. Let 
\begin{equation}
D_{N,l}:=\dim_{\Q}\gr^{\cL}_l\eZ^A_N.
\end{equation}
The computation of $D_{N,1}$ is relatively easy. We have
\begin{equation} \label{eqn:l1}
D_{N,1}=\begin{cases}
1, & \mbox{if $N>0$ is even}\\
0, & \mbox{else}.
\end{cases}
\end{equation}
More precisely, one can show (Corollary \ref{cor:diff}) that elliptic multiple zeta values of length one and weight $N$ are contained in $\Q \cdot (2\pi i)^N$ and vanish for odd $N$, and \eqref{eqn:l1} follows immediately from this. The computation of $D_{N,2}$, which is more elaborate, is the main result of this paper.
\begin{mthm}
We have
\begin{equation} \label{eqn:l2}
D_{N,2}=\begin{cases}
0, & \mbox{if $N$ is even}\\
\left\lfloor \frac{N}{3} \right\rfloor+1, & \mbox{if $N$ is odd}.
\end{cases}
\end{equation}
\end{mthm}
In the course of proving \eqref{eqn:l2}, we proceed as follows: Elliptic multiple zeta values satisfy two families of $\Q$-linear relations, namely the \textit{shuffle relations} and the \textit{Fay relations}~\cite{BMMS}. In order to count the number of linearly independent such relations, we introduce the \textit{Fay-shuffle space} $\FSh_2(N)$
\footnote{Independently, this space has already been studied by Brown in an unpublished preprint, in a slightly different context \cite{BrownAnatomy}.}
in length two, which is an elliptic analogue of the (depth two) double-shuffle space defined in~\cite{IKZ} in the context of multiple zeta values. By construction, the dimension of $\FSh_2(N)$ gives an upper bound for $D_{N,2}$, and using representation theory of the symmetric group $S_3$, we obtain for $\dim_{\Q}\FSh_2(N)$ precisely the numbers on the right hand side of \eqref{eqn:l2} (Theorem \ref{thm:Fayshuf}). In particular, the main theorem for even $N$ follows already from this.
In order to prove \eqref{eqn:l2} for odd $N$, we need a linear independence result for the elliptic double zeta values $I^A(r,N-r)$ for $r=0,\dots,k$, where $k=\lfloor \frac N3 \rfloor$ (Theorem \ref{thm:linind}). This uses Enriquez's differential equation for elliptic multiple zeta values and ends the proof of the main theorem. Moreover, since there are no non-trivial relations between elliptic double zeta values of different weights (Theorem \ref{thm:weight}), it also follows from the proof that every relation between elliptic double zeta values is a consequence of shuffle and Fay relations.

\bigskip

\noindent
1.5. {\bf Plan of the paper.}
In Section \ref{sec:prelim}, we give the definition of the differential forms $\omega^{(k)}$, following~\cite{BL}. Enriquez's definition of elliptic multiple zeta values is given in Section \ref{sec:eMZV}. We also indicate how his definition relates to homotopy invariant iterated integrals, and give an account of relations between elliptic multiple zeta values in the length two case \cite{BMS}.

In Section \ref{sec:dimeMZV}, we arrive at the main results of this paper. We introduce the algebra of A-elliptic multiple zeta values, and define the length filtration on it. The statement and proof of the main theorem \eqref{eqn:l2} then occupies the rest of Section \ref{sec:dimeMZV}. The appendix contains the computation of the determinant of a matrix with binomial coefficients, needed for a step in the proof of the main theorem, which we were unable to find in the literature.
\bigskip

\noindent
1.6. {\bf Acknowledgments.}
Very many thanks to Benjamin Enriquez for his very thorough reading of a first draft of this paper, as well as for many corrections and suggestions. Part of this paper was written during a research stay in Strasbourg in February 2015, and I would like to thank the IRMA Strasbourg for hospitality. Also, very many thanks to Francis Brown for many corrections and suggestions on an earlier version of this paper, for suggesting a simpler and more conceptual proof of Proposition \ref{prop:brown}, and for communicating to me his results on $\fu^{geom}$ \cite{BrownLetter}. Thanks also to Henrik Bachmann, Johannes Broedel, Ulf K\"uhn, Pierre Lochak, Oliver Schlotterer and Leila Schneps for helpful discussions and comments. This paper is part of the author's doctoral thesis at Universit\"at Hamburg, and I would like to thank my advisor Ulf K\"uhn for his constant support of my work.


\section{Preliminaries} \label{sec:prelim}

We begin by collecting some basic facts about iterated integrals, partly in order to fix our notation. Then we recall the construction of the doubly periodic differential forms $\omega^{(k)}$ on a once-punctured elliptic curve, which will be the building blocks of elliptic multiple zeta values. See also~\cite{BMMS} for an expository account.
\bigskip

\noindent
2.1. {\bf Iterated integrals.}
Let $\omega_1,\dots,\omega_r$ be a collection of smooth one-forms on a complex manifold $M$, and let $\gamma: [0,1] \rightarrow M$ be a piecewise smooth path. Write $f_i(t)\dd t$ for the pull back of $\omega_i$ along $\gamma$.
We define
\begin{equation}
\int_{\gamma}\omega_1\dots\omega_r:=\int_{0<t_1<\dots<t_r<1}f_1(t_1)\dd t_1\dots f_r(t_r)\dd t_r.
\end{equation}
and call this integral an \emph{iterated integral}. More generally, every $\C$-linear combination of iterated integrals as above will also be called an iterated integral. If $r=0$, we set $\int_{\gamma}=1$.


Iterated integrals satisfy some important algebraic relations, among them the shuffle product formula
\begin{equation} \label{eqn:shuf}
\int_{\gamma}\omega_1\dots\omega_r\int_{\gamma}\omega'_1\dots\omega'_s=\int_{\gamma}\omega_1\dots\omega_r \shuffle \omega'_1\dots\omega'_s,
\end{equation}
the composition of paths formula
\begin{equation} \label{eqn:comp}
\int_{\gamma_1\gamma_2} \omega_1\dots\omega_r=\sum_{k=0}^r\int_{\gamma_1}\omega_1\dots\omega_k \int_{\gamma_2}\omega_{k+1}\dots\omega_r,
\end{equation}
and the reversal of paths formula
\begin{equation} \label{eqn:revpath}
\int_{\gamma^{-1}}\omega_1\dots\omega_r=(-1)^r\int_{\gamma}\omega_r\dots\omega_1.
\end{equation}
For a more detailed account of iterated integrals, we refer to~\cite{Ch}. 
\bigskip

\noindent
2.2. {\bf The Kronecker series.}
Fix $\tau$ in the upper half-plane $\H:=\{\xi \in \C \, \vert \, \Im(\xi)>0\}$. We consider a version of the odd Jacobi theta function $\theta_{\tau}$, defined for $\xi \in \C$ by
\begin{equation}
\theta_{\tau}(\xi):=\sum_{n \in \Z}(-1)^nq^{\frac{1}{2}(n+\frac{1}{2})^2}e^{(n+\frac{1}{2})2\pi i\xi}, \quad q=e^{2\pi i\tau}.
\end{equation}
\begin{dfn}
The \textit{Kronecker series} $F_{\tau}: \C \times \C \rightarrow \C$~\cite{W,Z} is the meromorphic function defined by the formula
\begin{equation}
F_{\tau}(\xi,\alpha)=\frac{\theta_{\tau}'(0)\theta_{\tau}(\xi+\alpha)}{\theta_{\tau}(\xi)\theta_{\tau}(\alpha)}.
\end{equation}
\end{dfn}
The Kronecker series is quasi-periodic with respect to lattice translations in both variables
\begin{align}
F_{\tau}(\xi+1,\alpha)=F_{\tau}(\xi,\alpha), \quad F_{\tau}(\xi+\tau,w)=e^{-2\pi i\frac{\Im(\xi)}{\Im(\tau)}\alpha}F_{\tau}(\xi,\alpha), \\
F_{\tau}(\xi,\alpha+1)=F_{\tau}(\xi,\alpha), \quad F_{\tau}(\xi,\alpha+\tau)=e^{-2\pi i\frac{\Im(\alpha)}{\Im(\tau)}\xi}F_{\tau}(\xi,\alpha).
\end{align}
Viewed as a function of either variable, it has simple poles along the lattice $\Z+\Z\tau$. For more on the Kronecker series, see for example~\cite{W}, Part~II, and~\cite{Z}.
\bigskip

\noindent
2.3. {\bf Differential forms on an elliptic curve.}
Every complex elliptic curve admits a representation as a quotient of $\C$ by a lattice
\begin{equation}
E_{\tau}=\C/(\Z+\Z\tau),
\end{equation}
with $\tau \in \H$. We write $E^{\times}_{\tau}$ for $E_{\tau}$ with the point $0$ removed, and denote by $\xi=s+r\tau$, with $r,s \in \R$, the canonical coordinate on $E^{\times}_{\tau}$.

Consider the formal differential one-form (cf.~\cite{BL}, Section~3.5), where $\alpha$ is a formal variable
\begin{equation}
\Omega_{\tau}(\xi,\alpha):=e^{2\pi ir\alpha}F_{\tau}(\xi,\alpha)\dd \xi.
\end{equation}
\begin{prop} \label{prop:kro}
The one-form $\Omega_{\tau}(\xi,\alpha)$ satisfies the following properties:
\begin{enumerate}
\item
It is doubly periodic in $\xi$:
\begin{align}
\Omega_{\tau}(\xi+1,\alpha)=\Omega_{\tau}(\xi,\alpha), \quad \Omega_{\tau}(\xi+\tau,\alpha)=\Omega_{\tau}(\xi,\alpha).
\end{align}
Thus it descends to a formal differential one-form on $E^{\times}_{\tau}$.
\item
It is symmetric in its arguments, i.e.
\begin{equation}
\Omega_{\tau}(-\xi,-\alpha)=\Omega_{\tau}(\xi,\alpha).
\end{equation}
\item
It satisfies the Fay identity:
\begin{align}
\Omega_{\tau}(\xi_1,\alpha_1) \land \Omega_{\tau}(\xi_2,\alpha_2)&=\Omega_{\tau}(\xi_1-\xi_2,\alpha_1) \land \Omega_{\tau}(\xi,\alpha_1+\alpha_2) \notag \\
&+\Omega_{\tau}(\xi_2-\xi_1,\alpha_2) \land \Omega_{\tau}(\xi,\alpha_1+\alpha_2).
\end{align}
\item
We have the modular transformation formula
\begin{equation} \label{eqn:Omegamod}
\Omega_{\frac{a\tau+b}{c\tau+d}}((c\tau+d)^{-1}\xi,(c\tau+d)^{-1}\alpha)=(c\tau+d)\Omega_{\tau}(\xi,\alpha)
\end{equation}
for $\begin{pmatrix}a&b\\c&d\end{pmatrix} \in \SL_2(\Z)$.
\end{enumerate}
\end{prop}
\begin{prf}
Double periodicity and the Fay identity are deduced from analogous properties of the Kronecker function (cf.~\cite{BL}, Proposition~5). The symmetry is clear from the definition and from the fact that $\theta_{\tau}(\xi)$ is an odd function. Finally, the modular transformation formula follows from the modular transformation formula for $\theta_{\tau}(\xi)$ (cf.~\cite{Z}, Section~3, Theorem~(vi)).
\end{prf}
Following~\cite{BL}, Definition~7, we will consider $\Omega_{\tau}(\xi,\alpha)$ as a generating series of differential one-forms on $E^{\times}_{\tau}$. Recall that $F_{\tau}(\xi,\alpha)$ has a simple pole at $\alpha=0$, and that therefore $\Omega_{\tau}(\xi,\alpha)$ has a formal expansion in $\alpha$.
\begin{dfn}
Define a family $\{\omega^{(k)}\}_{k \geq 0}$ of real analytic differential one-forms on $E^{\times}_{\tau}$ by the formula
\begin{equation}
\Omega_{\tau}(\xi,\alpha)=\sum_{k \geq 0}\omega^{(k)}\alpha^{k-1}.
\end{equation}
\end{dfn}
Note that we have
\begin{equation}
\omega^{(k)}=f^{(k)}(\xi)\dd\xi,
\end{equation}
where each $f^{(k)}$ is a real analytic function on $E^{\times}_{\tau}$. The $\omega^{(k)}$ satisfy properties analogous to $\Omega_{\tau}(\xi,\alpha)$.
\begin{cor} \label{cor:fk}
The form $\omega^{(k)}$ has modular weight $k+1$, i.e.
\begin{equation}
\omega^{(k)}_{\vert \gamma}=(c\tau+d)^{k+1}\omega^{(k)}
\end{equation}
where the action of $\gamma \in \SL_2(\Z)$ on $\omega^{(k)}$ is defined as in~\eqref{eqn:Omegamod}. Moreover, $\omega^{(k)}=f^{(k)}(\xi)\dd\xi$ is doubly periodic and antisymmetric/symmetric depending on the weight
\begin{equation}
f^{(k)}(-\xi)=(-1)^kf^{(k)}(\xi),
\end{equation}
and we have the Fay identity
\begin{align}
f^{(m)}(\xi_1)f^{(n)}(\xi_2)&=-(-1)^nf^{(m+n)}(\xi_1-\xi_2) \notag\\
&+\sum_{r=0}^{n}\binom{m+r-1}{m-1}f^{(n-r)}(\xi_2-\xi_1)f^{(m+r)}(\xi_1) \notag\\
&+\sum_{r=0}^{m}\binom{n+r-1}{n-1}f^{(m-r)}(\xi_1-\xi_2)f^{(n+r)}(\xi_2).
\end{align}
\qed
\end{cor}

\bigskip

\noindent
2.4. {\bf Indefinite Eisenstein integrals.}
Recall the definition of the Eisenstein series $G_{2k}(\tau)$ for $k \geq 1$, where $\sigma_m(n)=\sum_{d | n}d^m$ is the $m$-th divisor sum:
\begin{equation}
G_{2k}(\tau)=2\zeta(2k)+2\frac{(2\pi i)^{2k}}{(2k-1)!}\sum_{n \geq 1}\sigma_{2k-1}(n)q^n, \quad q=e^{2\pi i\tau}.
\end{equation}
We extend this definition to $n=0$ by setting $G_0(\tau) \equiv -1$.
\begin{dfn} (\cite{MMV}, Example 4.10)
We define the \textit{indefinite Eisenstein integral} by
\begin{equation} \label{eqn:inteis}
\calg_{2k}(\tau)=\int_{\tau}^{i\infty}G_{2k}(\tau')-2\zeta(2k)d\tau'+\int_0^{\tau}2\zeta(2k)d\tau'.
\end{equation}
\end{dfn}
Note that since $G_{2k}(\tau)-2\zeta(2k) \sim O(e^{2\pi i\tau})$ as $\tau \to i\infty$, the left-hand integral in \eqref{eqn:inteis} is well-defined.

The following result is a direct consequence of $\C$-linear independence of the usual Eisenstein series. Denote by $\calo(\H)$ the $\C$-algebra of holomorphic functions on $\H$.
\begin{prop} \label{prop:gradeis}
Let $\Eis_{2\pi i} \subset \calo(\H)$ be the $\Q$-vector subspace spanned by all products
\begin{equation}
(2\pi i)^j\calg_{2k}(\tau), \quad j,k \geq 0,
\end{equation}
and define a grading on $\Eis_{2\pi i}$ by giving $(2\pi i)^j\calg_{2k}(\tau)$ weight $j+2k$. Then $\Eis_{2\pi i}$ is a graded $\Q$-vector space.
\end{prop}

\section{Review of elliptic multiple zeta values} \label{sec:eMZV}

In this section, we introduce Enriquez's elliptic multiple zeta values and recall some of their basic properties. The important differential equation, which elliptic multiple zeta values satisfy is also stated. This differential equation, due to Enriquez (cf.~\cite{E2}, and also~\cite{CEE,LR}), shows in particular that elliptic multiple zeta values of length one are constant, and that elliptic multiple zeta values of length two are linear combinations of indefinite integrals of Eisenstein series \eqref{eqn:inteis}. We end this section by giving some results on $\Q$-linear relations satisfied by elliptic multiple zeta values in the special case of length two. We refer to~\cite{BMMS,BMS} for a study of more general relations between elliptic multiple zeta values in higher lengths.
\bigskip

\noindent
3.1. {\bf Definition and first properties.}
Let $\tau \in \H$, and consider the complex elliptic curve $E_{\tau}$. Its universal covering map is given by the projection $\C \rightarrow \C/\Z+\Z\tau \cong E_{\tau}$. We define paths $\alpha,\beta$ on $E_{\tau}$ to be the images of the straight line paths from $0$ to $1$ and from $0$ to $\tau$ respectively.

Elliptic multiple zeta values are defined as iterated integrals of the formal differential one-form 
\begin{equation}
\ad(x)\Omega_{\tau}(\xi,\ad(x))(y)
\end{equation}
defined on the once-punctured elliptic curve $E^{\times}_{\tau}:=E_{\tau} \setminus \{0\}$ with values in the free Lie algebra on two generators $x,y$, along the paths $\alpha$ and $\beta$ (cf.~\cite{E}, D\'efinition~2.6, eq.~(18), and~\cite{BMMS}, Section~2.2.1). By~\cite{BL}, eq.~(3.8), the above differential form has a simple pole at $0 \in E_{\tau}$ with residue $2\pi i\ad(x)(y)$, in the sense that
\begin{equation} \label{eqn:residue}
\lim_{\xi \to 0}\xi\ad(x)\Omega_{\tau}(\xi,\ad(x))(y)=2\pi i\ad(x)(y),
\end{equation}
This can be done choosing non-zero tangent vectors at $0$, i.e. tangential base points in the sense of~\cite{Del}, \S~15, and leads to the following definition.
\begin{dfn} \label{dfn:eMZV}
We define the A-elliptic multiple zeta value $I^A(n_1,\dots,n_r;\tau)$\footnote{In \cite{BMMS}, this is denoted by $\omega(n_1,...,n_r)$} to be the coefficient of $\ad^{n_1}(x)(y)\dots\ad^{n_r}(x)(y)$ in the regularized generating series
\begin{equation}
\lim_{t \to 0}(-2\pi it)^{\ad(x)(y)}\exp\bigg[\int_{\alpha_t^{1-t}}\ad(x)\Omega_{\tau}(\xi,\ad(x)(y))\bigg](-2\pi it)^{-\ad(x)(y)},
\end{equation}
where $\alpha_t^{1-t}$ denotes the restriction of the path $\alpha$ to the interval $[t,1-t]$, and for a differential one-form $\omega$, we write
\begin{equation}
\exp\bigg[\int \omega\bigg]=1+\sum_{k \geq 1}\int \underbrace{\omega\dots\omega}_{\textrm{k-times}}.
\end{equation}
The regularization of the integral is performed as in~\cite{Del}, Proposition 15.45 (with the order of integration reversed), with respect to the tangent vector $-\frac{\partial}{\partial z}=(-2\pi i)^{-1}\frac{\partial}{\partial \xi}$ at $0$ at the starting point of the path of integration (where $z=e^{2\pi i\xi}$), and with respect to the tangent vector $\frac{\partial}{\partial z}=(2\pi i)^{-1}\frac{\partial}{\partial \xi}$ at $0$ at the endpoint of the path of integration.
\end{dfn}
$I^A(n_1,\dots,n_r;\tau)$ is said to have \textit{length} $r$ and \textit{weight} $n_1+\dots+n_r$. If $n_1,n_r \neq 1$, then we simply have
\begin{equation}
I^A(n_1,\dots,n_r;\tau)=\int_{\alpha} \omega^{(n_1)}\dots\omega^{(n_r)},
\end{equation}
where the forms $\omega^{(k)}$ have been defined in Section \ref{sec:prelim}. Note that this is well-defined, since $\omega^{(k)}$ does not have a pole at the puncture $0$, if $k \neq 1$, by \eqref{eqn:residue}.

The following result has essentially already been proved in~\cite{E}, Proposition~2.8. Since our conventions differ slightly from \cite{E}, we repeat the proof for completeness.
\begin{prop} \label{prop:compass}
For all $n_1,\dots,n_r \geq 0$, we have that $I^A(n_1,\dots,n_r;\tau)$ is equal to $(-1)^r$ times the coefficient of the word $\ad^{n_r}(x)(y)...\ad^{n_1}(x)(y)$ in the series
\begin{equation}
e^{-\pi i\ad(x)(y)}A(\tau),
\end{equation}
where $A(\tau)$ is Enriquez's A-associator (cf.~\cite{E2}, Section~5.2).
\end{prop}
\begin{prf}
Recall that $A(\tau)$ is defined as
\begin{equation}
A(\tau)=X(z)^{-1}X(z+1),
\end{equation}
where $X(z)$ is the unique holomorphic function defined on $\{z \in \C \, | \, z=u+vi, \ u \mbox{ or } v \in (0,1)\} \subset \C$, which solves the equation 
\begin{equation}
\dd X(z)=-\ad(x)\Omega_{\tau}(\xi,\ad(x))(y) \cdot X(z),
\end{equation}
and satisfies $X(z) \sim (-2\pi iz)^{-\ad(x)(y)}$, as $z \to 0$, where we choose the principal branch of the logarithm so that $\log(\pm i)=\pm \frac{\pi i}{2}$. Using Picard iteration, and standard properties of iterated integrals (cf.~Section~2.1), we see that $A(\tau)$ equals
\begin{align}
\bigg(\lim_{t \to 0}(-2\pi it)^{\ad(x)(y)}\exp\bigg[\int_{\alpha_t^{1-t}}(-\ad(x)\Omega_{\tau}(\xi,\ad(x)(y)))\bigg]e^{\pi i\ad(x)(y)}(-2\pi it)^{-\ad(x)(y)}\bigg)^{op},
\end{align}
where the superscript $op$ denotes the opposite multiplication on the algebra $\C \langle\langle x,y\rangle\rangle$, defined by $(f \cdot g)^{op}=g \cdot f$. Now multiplying by $e^{-\pi i\ad(x)(y)}$, the result follows.
\end{prf}
Proposition \ref{prop:compass} can be seen as an analogue for A-elliptic multiple zeta values of the fact that multiple zeta values are the coefficients of the Drinfel'd associator~\cite{Dr,LM}.
\begin{rmk}
In~\cite{E}, one also finds a notion of B-elliptic multiple zeta value, which we denote here by $I^B(n_1,\dots,n_r;\tau)$, defined as iterated integrals along the other homology cycle $\beta$ of an elliptic curve. Similar to A-elliptic multiple zeta values, B-elliptic multiple zeta values can be obtained as coefficients of Enriquez's B-associator~\cite{E2}. The modular transformation formula (cf.~\cite{E}, eq.~(26))
\begin{equation}
I^B(n_1,\dots,n_r;\tau)=\tau^{-(n_1+\dots+n_r)+r}I^A(n_1,\dots,n_r;-\tau^{-1})
\end{equation}
establishes a close connection between A-elliptic and B-elliptic multiple zeta values. A more detailed treatment of B-elliptic multiple zeta values will be given elsewhere.
\end{rmk}
By Proposition~5.3 of~\cite{E}, every A-elliptic multiple zeta value has a Fourier expansion
\begin{equation}
\sum_{k \geq 0}a_kq^k, \quad q=e^{2\pi i\tau},
\end{equation}
with the coefficients $a_k$ contained in $\cZ[(2\pi i)^{-1}]$, where $\cZ$ denotes the $\Q$-algebra of multiple zeta values. From the asymptotic behavior of the KZB associator (cf.~\cite{BMS}, Section~2.3), one can read off the constant term in the above $q$-expansion. We will do this explicitly in lengths one and two.
\begin{prop} \label{prop:constterm}
Denote by $c_{n_1,\dots,n_r}$ the constant term in the Fourier expansion of $I^A(n_1,\dots,n_r;\tau)$. It is
\begin{equation}
c_n=
\begin{cases}
0 ,& \mbox{for $n=1$}\\
\frac{B_n}{n!}(2\pi i)^n, & \mbox{else},
\end{cases}
\end{equation}
and
\begin{equation}
c_{m,n}=
\begin{cases}
0, & \mbox{if $m=n=1$}\\
-\frac 12 (2\pi i)^{m+1}\frac{B_mB_1}{m!},  & \mbox{if $n=1$ and $m \neq 1$}\\
\frac 12 (2\pi i)^{m+n}\frac{B_mB_n}{m!n!} ,& \mbox{else},
\end{cases}
\end{equation}
where $B_n$ are the Bernoulli numbers defined by the series $\sum_{n \geq 0}\frac{B_n}{n!}t^n=\frac{t}{e^t-1}$.
\end{prop}
\begin{prf}
By~\cite{E} (cf.~also~\cite{BMS}, eq.~(2.43)), we have
\begin{equation} \label{eqn:constterm}
\lim_{\tau \to i\infty} e^{\pi it}A(\tau)=e^{\pi it}\Phi(\tilde{y},t)e^{2\pi i\tilde{y}}\Phi(\tilde{y},t)^{-1},
\end{equation}
where $t=-\ad(x)(y)$ and $\tilde{y}=-\frac{\ad(x)}{e^{2\pi i\ad(x)}-1}(y)$. Now by Proposition \ref{prop:compass}, we see that $c_{n_1,\dots,n_r}$ is given by the coefficient of $(-1)^r\ad^{n_r}(x)(y)\dots\ad^{n_1}(x)(y)$ on the left hand side of \eqref{eqn:constterm}, and the result follows by comparing coefficients.
\end{prf}
\bigskip

\noindent
3.2. {\bf Homotopy invariance.}
Although not strictly needed for the rest of this paper, we show in this section that the A-elliptic multiple zeta values, which were defined as iterated integrals on $E^{\times}_{\tau}$ are in fact given by \textit{homotopy invariant} iterated integrals. This is non-trivial, since the differential forms $\omega^{(k)}$ are not closed for $k \geq 1$, and thus already the value of the ordinary line integral
\begin{equation}
\int_{\gamma}\omega^{(k)}, \quad k \geq 1
\end{equation}
is not invariant under homotopies of the path $\gamma$.

The fundamental result we will need is the following. Recall that $\xi=s+r\tau$ with $r,s \in \R$ is the canonical coordinate on $E^{\times}_{\tau}=\C/(\Z+\Z\tau) \setminus \{0\}$.
\begin{thm}[Brown-Levin] \label{thm:BL}
Every homotopy invariant iterated integral on $E^{\times}_{\tau}$ is a $\C$-linear combination of coefficients of the words $x,y$ in the generating series
\begin{equation} \label{eqn:geniter}
T=1+\sum_{k \geq 1}\int \underbrace{J\dots J}_{\textrm{k-times}},
\end{equation}
where 
\begin{equation} \label{eqn:kzbeqn}
J=-\nu x+\ad(x)\Omega_{\tau}(\xi,\ad(x))(y)=-\nu x+\sum_{k \geq 0}\omega^{(k)}\ad^k(x)(y)
\end{equation} 
with $\nu=2\pi i \dd r$.
\end{thm}
Since every coefficient of $T$ is a homotopy invariant iterated integral, for any pair of points $\xi,\rho \in E^{\times}_{\tau}$, we can consider $T$ as a function
\begin{align}
T(\_ ): \pi_1(E^{\times}_{\tau};\xi,\rho) &\rightarrow \C\langle\langle x,y\rangle\rangle\\
\gamma &\mapsto T(\gamma),
\end{align}
where $T(\gamma)$ denotes integration of $T$ along the path $\gamma$.
Using Deligne's tangential base points (cf.~\cite{Del}, \S~15), we can extend the domain of $T(\_ )$ to also include paths which start or end at $0$. In the case of the path $\alpha$ (cf.~Section~3.1) and the tangent vectors $(-2\pi i)^{-1}$ and $(2\pi i)^{-1}$ at $0$, this gives
\begin{equation}
T^{\Reg}(\alpha)=\lim_{t \to 0}(-2\pi i)^{\ad(x)(y)}\big[T(\alpha^{1-t}_t)\big](-2\pi i)^{-\ad(x)(y)}.
\end{equation}
We now compare with the A-elliptic multiple zeta values defined in the last section. Since the differential form $\nu$ vanishes along the path $\alpha$, we see from Theorem \ref{thm:BL} that every coefficient $T^{\Reg}(\alpha)_w$ of $T^{\Reg}(J;\alpha)$ is given by a special $\Z$-linear combination of A-elliptic multiple zeta values $I^A(n_1,\dots,n_r;\tau)$.

For example, if $w$ is the word $yx^{n_1}\dots yx^{n_r}$ and $\underline{n}=(n_1,\dots,n_r)$, then the coefficient $T^{\Reg}(\alpha)_w$ is obtained from expanding Lie monomials $\ad^k(x)(y)$ into linear combinations of words in $x,y$, and then collecting the respective terms in \eqref{eqn:kzbeqn}. The result is
\begin{equation} \label{eqn:transl}
T^{\Reg}(\alpha)_w=(-1)^{n_1+\dots+n_r}I^A(\underline{n};\tau)+\sum_{\underline{n'} \prec \underline{n}} c_{\underline{n'}}I^A(\underline{n'};\tau),
\end{equation}
for integers $c_{\underline{n'}} \in \Z$ coming from the expansion of Lie monomials $\ad^{k}(x)(y)$. Here and in the following, the multi-indices $\underline{n} \in \N^r$ being ordered lexicographically for the natural order on $\N$.
\begin{prop} \label{prop:transl}
For every collection of integers $n_1,\dots,n_r \geq 0$, the A-elliptic multiple zeta value $I^A(n_1,\dots,n_r;\tau)$ is a $\Z$-linear combination of coefficients of $T^{\Reg}(\alpha)_w$.
\end{prop}
In particular, this proves that every A-elliptic multiple zeta value is the restriction of a homotopy invariant iterated integral.
\begin{prf}
We use induction on the set of multi-indices of a fixed length $r$, ordered lexicographically. For $\underline{n}=(0,\dots,0) \in \N^r$, it is clear that
\begin{equation}
I^A(0,\dots,0;\tau)=T^{\Reg}(\alpha)_{y^r}.
\end{equation}
Now let $\underline{n}=(n_1,\dots,n_r)$ be a multi-index of length $r$ and assume the proposition for all multi-indices $\underline{n'}$ such that $\underline{n'} \prec \underline{n}$ in the lexicographical ordering. Now by \eqref{eqn:transl}, we have
\begin{equation}
I^A(\underline{n};\tau)=(-1)^{n_1+\dots+n_r}\left(T^{\Reg}(\alpha)_w-\sum_{\underline{n'} \prec \underline{n}} c_{\underline{n'}}I^A(\underline{n'};\tau) \right),
\end{equation}
where $w=yx^{n_1}\dots yx^{n_r}$. The terms on the right hand side of the last equation can be expressed as coefficients of $T^{\Reg}(\alpha)$, by the induction hypothesis, and this ends the proof.
\end{prf}
\bigskip

\noindent
3.3. {\bf The differential equation for A-elliptic multiple zeta values.}
Write
\begin{equation}
\cali^A(X_1,\dots,X_r;\tau):=\sum_{n_1,\dots,n_r \geq 0}I^A(n_1,\dots,n_r;\tau)X^{n_1-1}_1\dots X^{n_r-1}_r
\end{equation}
for the generating series of A-elliptic multiple zeta values of length $r$. The derivative of $\cali^A$ is computed in~\cite{E}, Th\'eor\`eme~3.10. Precisely:
\begin{thm}[Enriquez] \label{thm:diff}
For all $r \geq 0$, we have
\begin{align} \label{eqn:diffenr}
2\pi i\frac{\partial}{\partial \tau}\cali^A(X_1,\dots,X_r;\tau)&=\wp_{\tau}^*(X_1)\cali^A(X_2,\dots,X_r;\tau)-\wp_{\tau}^*(X_r)\cali^A(X_1,\dots,X_{r-1};\tau) \notag\\
&+\sum_{i=1}^{r-1}(\wp_{\tau}^*(X_{i+1})-\wp_{\tau}^*(X_i))\cali^A(X_1,\dots,X_i+X_{i+1},\dots,X_r;\tau),
\end{align}
where $\wp_{\tau}^*(\alpha)=\sum_{k=-1}^{\infty}(2k+1)G_{2k+2}(\tau)\alpha^{2k}$ considered as a formal power series in $\alpha$.

\qed
\end{thm}
By comparing coefficients on both sides of \eqref{eqn:diffenr}, one can extract explicit formulas for the $\tau$-derivatives of individual A-elliptic multiple zeta values $I^A(n_1,\dots,n_r;\tau)$. Together with Proposition \ref{prop:constterm}, this leads to expressions for A-elliptic multiple zeta values as iterated integrals of Eisenstein series (cf.~\cite{BMS}, Section~4.2). We only do the cases of length one and two here.
\begin{cor} \label{cor:diff}
Elliptic multiple zeta values of length one are constant, i.e.
\begin{equation} \label{eqn:cor1}
2\pi i\frac{\partial}{\partial \tau}I^A(n;\tau)=0.
\end{equation}
More precisely, writing $I^A(n)$ instead of $I^A(n;\tau)$, we have $I^A(1)=0$ and
\begin{equation} \label{eqn:cor1b}
I^A(n)=\frac{(2\pi i)^nB_n}{n!}, \quad n \neq 1.
\end{equation}
In length two, we have
\begin{align} \label{eqn:cor2a}
2\pi i\frac{\partial}{\partial \tau}I^A(0,n;\tau)&=(-1)^n2\pi i\frac{\partial}{\partial \tau}I^A(n,0;\tau) \notag \\
&=-nG_{n+1}(\tau)I^A(0)+nG_0(\tau)I^A(n+1),
\end{align}
and if $m,n \neq 0$
\begin{align} \label{eqn:cor2b}
2\pi i\frac{\partial}{\partial \tau}I^A(m,n;\tau)&=-nG_{n+1}(\tau)I^A(m)+mG_{m+1}(\tau)I^A(n) \notag \\
&-(-1)^m(m+n)G_{m+n+1}(\tau)I^A(0) \notag \\
&+\sum_{k=1}^{m+n+1}(m+n-k)\left( \binom{k-1}{m-1}-\binom{k-1}{n-1} \right)G_{m+n+1-k}(\tau)I^A(k),
\end{align}
where we recall that by definition $G_n(\tau) \equiv 0$ whenever $n$ is odd. In particular, $\frac{\partial}{\partial \tau}I^A(m,n;\tau)=0$ if $m+n$ is even, and is given by
\begin{equation} \label{eqn:cor2c}
I^A(m,n;\tau)=\frac{(2\pi i)^{m+n}}{2} \frac{B_mB_n}{m!n!}.
\end{equation}
\end{cor}
\begin{prf}
\eqref{eqn:cor1} and~\eqref{eqn:cor1b} follow directly from Theorem \ref{thm:diff} and Proposition \ref{prop:constterm} respectively. Similarly,~\eqref{eqn:cor2a} and~\eqref{eqn:cor2b} follow from Theorem \ref{thm:diff} by comparing coefficients, and~\eqref{eqn:cor2c} follows from Proposition \ref{prop:constterm} bearing in mind that $I^A(2j+1)=0$ for all $j \geq 0$ by~\eqref{eqn:cor1b}.
\end{prf}

\bigskip

\noindent
3.4. {\bf Relations in length two.}
In this section, we review several $\Q$-linear relations between A-elliptic multiple zeta values, restricting ourselves to the length two case. More general formulas can be found in~\cite{BMMS,BMS}.
\begin{prop} \label{prop:Qrel}
For all $m,n \geq 0$, we have the following relations between elliptic double zeta values:
\begin{enumerate}
\item (Reflection relation)
\begin{equation}
I^A(m,n;\tau)=(-1)^{m+n}I^A(n,m;\tau)
\end{equation}
\item (Shuffle relation)
\begin{equation}
I^A(m,n;\tau)+I^A(n,m;\tau)=I^A(m)I^A(n)
\end{equation}
\item (Fay relation)
\begin{align}
I^A(m,n;\tau)&=-\delta_{m,1}\delta_{n,1}3\zeta(2)-(-1)^nI^A(0,m+n;\tau) \notag\\
&+\sum_{r=0}^n(-1)^{n-r}\binom{m-1+r}{m-1}I^A(m+r,n-r;\tau) \notag \\
&+\sum_{r=0}^m(-1)^{n+r}\binom{n-1+r}{n-1}I^A(n+r,m-r;\tau)
\end{align}
\end{enumerate}
\end{prop}
Note that the right hand side of ii) is contained in $\cL_1 \eZ^A$ (cf.~Corollary~\ref{cor:diff}(i)), and therefore vanishes in $\gr^{\cL}_2\eZ^A$. Therefore, ii) can be seen as a $\Q$-linear relation in $\gr^{\cL}_2\eZ^A$.
\begin{prf}
\begin{enumerate}
\item
This follows from the reversal of paths formula \eqref{eqn:revpath}, using in addition the parity properties of the functions $f^{(n)}$ (Corollary \ref{cor:fk}).
\item
This is just a special case of the usual shuffle product property of iterated integrals \eqref{eqn:shuf}.
\item
If $m=n=1$, then we get $I^A(1,1;\tau)=-2I^A(1,1;\tau)+I^A(0,2;\tau)+2I^A(2,0;\tau)-3\zeta(2)$, which holds by Corollary \ref{cor:diff}.
Otherwise, by i) we can assume that $n \neq 1$.

Choose $\varepsilon>0$ and consider the function
\begin{equation}
\Xi^{m,n}_{\varepsilon}(x)=\int_{\varepsilon}^{x}f^{(n)}(\xi_2-x)\int_{\varepsilon}^{\xi_2}f^{(m)}(\xi_1)\dd\xi_1\dd\xi_2,
\end{equation}
for $x \in [\varepsilon,1]$. Since $n \neq 1$, $f^{(n)}$ is smooth on all of $[0,1]$, and thus $\Xi^{m,n}_{\varepsilon}$ is smooth as well. Moreover $\lim_{\varepsilon \to 0} \Xi^{m,n}_{\varepsilon}(1)=I^A(m,n;\tau)$, due to the periodicity of $f^{(n)}$ (Corollary \ref{cor:fk}). Now
\begin{equation}
\Xi^{m,n}_{\varepsilon}(t)=\int_{\varepsilon}^t(\Xi^{m,n}_{\varepsilon})'(x)\dd x=\int_{\varepsilon}^t \int_{\varepsilon}^{x}f^{(n)}(\xi_2-x)f^{(m)}(\xi_2)\dd\xi_2\dd x.
\end{equation}
Using the Fay identity for the $f^{(k)}$ (Corollary \ref{cor:fk}) we get
\begin{align}
\Xi^{m,n}_{\varepsilon}(t)=\int_{\varepsilon}^t \int_{\varepsilon}^{x}&\left\{ -(-1)^nf^{(m+n)}(x)+\sum_{r=0}^n\binom{m-1+r}{m-1}f^{(n-r)}(-x)f^{(m+r)}(\xi_2) \right. \notag \\
& \left.+\sum_{r=0}^m\binom{n-1+r}{n-1}f^{(n+r)}(\xi_2-x)f^{(m-r)}(x) \right\}\dd\xi_2\dd x.
\end{align}
Now we evaluate both sides at $t=1$ and pass to the limit $\varepsilon \to 0$ to obtain the result.
\end{enumerate}
\end{prf}
It will be important for the sequel to write out the relations in Proposition \ref{prop:Qrel} using generating series. This gives (cf.~\cite{BMS}, Section~2.2):
\begin{enumerate}
\item \textit{(Reflection relation)}
\begin{equation}
\cali^A(X,Y;\tau)=\cali^A(-Y,-X;\tau)
\end{equation}
\item \textit{(Shuffle relation)}
\begin{equation}
\cali^A(X,Y;\tau)+\cali^A(Y,X;\tau)=\cali^A(X;\tau)\cali^A(Y;\tau)
\end{equation}
\item \textit{(Fay relation)}
\begin{equation}
\cali^A(X,Y;\tau)+\cali^A(X+Y,-Y;\tau)+\cali^A(-X-Y,X;\tau)=-3\zeta(2)
\end{equation}
\end{enumerate}
\section{The dimension of the space of elliptic double zeta values} \label{sec:dimeMZV}

We define the main algebraic object of this paper, namely the $\Q$-vector space $\eZ^A$ of A-elliptic multiple zeta values. A proof parallel to the one for normal multiple zeta values, using the shuffle product, reveals that $\eZ^A$ is in fact a $\Q$-algebra, which is naturally filtered by the length of A-elliptic multiple zeta values. We then state our main result about the dimensions of the weight $N$ component of $\eZ^A$ in length two, namely Theorem \ref{mthm}. The rest of the section is then devoted to proving this theorem.
\bigskip

\noindent
4.1. {\bf The algebra of A-elliptic multiple zeta values}

\begin{dfn}
Define the $\Q$-vector space of A-elliptic multiple zeta values to be
\begin{equation}
\eZ^A=\langle I^A(n_1,\dots,n_r;\tau) \, \vert \, n_1,\dots,n_r \geq 0 \rangle_{\Q} \subset \calo(\H).
\end{equation}
Also, we define its weight $N$ component to be
\begin{equation}
\eZ^A_N=\langle I^A(n_1,\dots,n_r;\tau) \, \vert \, n_1+\dots +n_r=N \rangle_{\Q} \subset \calo(\H).
\end{equation}
\end{dfn}

\begin{prop}
The $\Q$-vector space $\eZ^A$ is a $\Q$-subalgebra of $\calo(\H)$.
\end{prop}
\begin{prf}
By definition, $\eZ^A$ is linearly spanned by the coefficients of the series
\begin{equation} \label{eqn:genseries}
\lim_{t \to 0}(-2\pi it)^{\ad(x)(y)}\exp\bigg[\int_{\alpha_t^{1-t}}\ad(x)\Omega_{\tau}(\xi,\ad(x)(y))\bigg](-2\pi it)^{-\ad(x)(y)}.
\end{equation}
By \eqref{eqn:comp}, the series \eqref{eqn:genseries} is group-like for the unique coproduct on $\C \langle \langle x,y \rangle \rangle$ for which $x$ and $y$ are primitive, which is equivalent to its coefficients satisfying the shuffle product formula. But since the elements $\ad^n(x)(y)$ are also primitive, the proposition follows.
\end{prf}
\bigskip

\noindent
4.2. {\bf The length filtration.}
\begin{dfn}
We denote by
\begin{equation}
\cL_l\eZ^A_N=\langle I^A(n_1,\dots,n_r;\tau) \in \eZ^A_N \, \vert \, r \leq l \rangle_{\Q}
\end{equation}
the $\Q$-subspace of A-elliptic multiple zeta values of weight $N$ and of length at most $l$. We also denote by
\begin{equation}
\gr^{\cL}_l\eZ^A_N
\end{equation}
the associated graded.
\end{dfn}
It is straightforward to extend the length filtration to all of $\eZ^A$. Also note that the length of $I^A(n_1,\dots,n_r;\tau)$ is just its length as an iterated integral (cf.~Definition~\ref{dfn:eMZV}).
\begin{prop}
$\eZ^A$ is a filtered $\Q$-algebra for the length filtration. More precisely, the multiplication law on $\eZ^A$ satisfies
\begin{equation}
\cL_l\eZ^A_N \cdot \cL_{l'}\eZ^A_{N'} \subset \cL_{l+l'}\eZ^A_{N+N'},
\end{equation}
for all $l,l',N,N' \geq 0$.
\end{prop}
\begin{prf}
This follows directly from the shuffle product formula \eqref{eqn:shuf} and the definition of the weight for elliptic multiple zeta values.
\end{prf}
In analogy with the weight grading conjecture for multiple zeta values, due to Zagier, a natural conjecture is the following.
\begin{conj} \label{conj:weight}
$\eZ^A$ is graded for the weight, i.e. the natural morphism
\begin{align}
\bigoplus_{N \geq 0} \eZ^A_N &\rightarrow \calo(\H)\\
(I^A_N)_{N \geq 0} &\mapsto \sum_{N \geq 0} I^A_N,
\end{align}
with $I^A_N \in \eZ^A_N$, induced by the inclusions $\eZ^A_N \subset \calo(\H)$, is injective.
\end{conj}

\bigskip

\noindent
4.3. {\bf Statement of the main results.}
The next theorem is the first main result of this paper, and provides a partial resolution of Conjecture \ref{conj:weight}.
\begin{thm} \label{thm:weight}
Denote by
$
\bigoplus_{N \geq 0} \cL_2\eZ^A_N
$
the outer direct sum of the $\Q$-vector spaces $\cL_2\eZ^A_N$. The natural map
\begin{equation}
\bigoplus_{N \geq 0} \cL_2\eZ^A_N \rightarrow \calo(\H)
\end{equation}
induced by the inclusions $\cL_2\eZ^A_N \subset \calo(\H)$ is injective.
\end{thm}
In other words, every $\Q$-linear relation between elliptic double zeta values can be decomposed into a family of $\Q$-linear relations, one for each weight $N$.
\begin{prf}
By Proposition \ref{prop:constterm} and Corollary \ref{cor:diff}, we have
\begin{equation}
\cL_2\eZ^A_N \subset \Q[2\pi i]_N \oplus (\Eis_{2\pi i})_{N+1},
\end{equation}
where a subscript $N$ denotes the weight $N$ component.
The theorem now follows from the fact that the $\Q$-vector subspace $\bigoplus_{N \geq 0}\Q[2\pi i]_N \oplus (\Eis_{2\pi i})_{N+1} \subset \calo(\H)$ is graded by Proposition \ref{prop:gradeis}, using the transcendence of $\pi$.
\end{prf}
We state our second main result.
\begin{thm} \label{mthm}
\begin{enumerate}
\item
Let $N \geq 0$ and $D_{N,2}:=\gr^{\cL}_2\eZ^A_N$. Then
\begin{equation}
D_{N,2}=\begin{cases}
0, & \mbox{if $N$ is even}\\
\left \lfloor \frac N3 \right \rfloor+1, & \mbox{if $N$ is odd}.
\end{cases}
\end{equation}
\item
Every $\Q$-linear relation between elliptic double zeta values is a consequence of Fay and shuffle relations.
\end{enumerate}
\end{thm}
The proof of Theorem \ref{mthm} will occupy the rest of this paper, and is completed in the last section. The idea is to count relations between elliptic double zeta values, which gives an upper bound for the dimension of $\gr^{\cL}_2\eZ^A_N$. The second step is then to show that this upper bound is attained by the elliptic double zeta values, which is achieved by using the explicit expression of elliptic double zeta values as integrals of Eisenstein series.
\bigskip

\noindent
4.4. {\bf The Fay-shuffle space.}
In this section, we introduce and study the Fay-shuffle space. Independently, this space has already been defined by Brown in the context of constructing polar solutions to the linearized double shuffle equations \cite{BrownAnatomy}.

Denote by $V'_N \subset \Q(X,Y)_N$ the subspace of rational functions, which are homogeneous of degree $N-2$, and have at most simple poles along $X=0$ and $Y=0$, and no other poles:
\begin{equation}
V'_N=\{f \in \Q(X,Y)_{N-2} \, \vert \, XY \cdot f \in \Q[X,Y]\}.
\end{equation}
\begin{dfn}
For $N \geq 0$, we define the length two \emph{Fay-shuffle space}, $\FSh_2(N)$ as the set of rational functions $P \in V'_N$, which satisfy
\begin{align} \label{eqn:fsh}
P(X,Y)+P(X+Y,-Y)+P(-X-Y,X)&=0, \quad P(X,Y)+P(Y,X)=0.
\end{align}
\end{dfn}
By substituting $X \mapsto -Y$ and $Y \mapsto -X$, one sees that $P \in \FSh_2(N)$ necessarily satisfies the reflection relation $P(X,Y)=P(-Y,-X)$. This immediately implies the
\begin{prop} \label{prop:zeroeven}
If $N$ is even, then $\FSh_2(N)=\{0\}$.
\end{prop}
\begin{prf}
Every $P \in \FSh_2(N)$ satisfies
\begin{equation} \label{eqn:parity}
P(X,Y)=P(-Y,-X), \quad P(X,Y)=-P(Y,X)
\end{equation}
hence $P \equiv 0$, if $N$ is even.
\end{prf}
The significance of the Fay-shuffle space is that its dimension gives an upper bound for the dimension of $\gr^{\cL}_2\eZ^A_N$.
\begin{prop} \label{rmk:Fayshuf}
For all $N \geq 0$, We have
\begin{equation}
\dim_{\Q} \gr^{\cL}_2\eZ^A_N \leq \dim_{\Q} \FSh_2(N).
\end{equation}
\end{prop}
\begin{prf}
Let $(V'_N)^*$ be the dual space of $V'_N$. Then, since the elements of $\gr^{\cL}_2\eZ^A_N$ satisfy the defining equations of $\FSh_2(N)$ (cf.~Proposition~\ref{prop:Qrel}; note that $-3\zeta(2)=\frac{3}{2}I^A(2;\tau)$ by Corollary \ref{cor:diff} and thus vanishes in $\gr^{\cL}_2\eZ^A_2$), the natural surjection
\begin{align}
(V'_N)^* &\rightarrow \gr^{\cL}_2\eZ^A_N \notag \\
(X^{k-1}Y^{l-1})^* &\mapsto I^A(k,l;\tau) \mod \cL_1\eZ^A_N
\end{align}
factors through the annihilator $\FSh^0_2(N) \subset (V'_N)^*$ of $\FSh_2(N)$, and therefore
\begin{equation}
\dim_{\Q} \gr^{\cL}_2\eZ^A_N \leq \dim_{\Q} \left[(V'_N)^*/(\FSh^0_2(N))\right]=\dim_{\Q} \FSh_2(N).
\end{equation}
\end{prf}
\bigskip

\noindent
4.5. {\bf Computation of the Fay-shuffle space.}
The goal of this section is to prove the following
\begin{thm} \label{thm:Fayshuf}
We have
\begin{equation}
\dim_{\Q}\FSh_2(N)=
\begin{cases}
0, & \mbox{if $N$ is even}\\
\left\lfloor \frac N3 \right\rfloor+1, & \mbox{if $N$ is odd.}
\end{cases}
\end{equation}
\end{thm}
For even $N$, the theorem follows from Proposition \ref{prop:zeroeven}. For odd $N$, the proof is divided into two propositions. First, if $N$ is odd, the space $\FSh_2(N)$ splits into a polynomial part, and a non-polynomial part as follows. Let $V_N \subset \Q[X,Y]$ be the subspace of homogeneous polynomials of degree $N$.
\begin{prop} \label{prop:decomp}
For odd $N \geq 1$, we have
\begin{equation}
\FSh_2(N) \cong \FSh_2(N)^{pol} \oplus \Q \widetilde{P},
\end{equation}
where 
\begin{equation}
\widetilde{P}(X,Y)=\frac{X^{N-1}}{Y}-\frac{Y^{N-1}}{X}-\frac{X^{N-1}-Y^{N-1}}{X+Y},
\end{equation}
and $\FSh_2(N)^{pol}=\FSh_2(N) \cap V_{N-2}$ denotes the polynomial part of the Fay-shuffle space.
\end{prop}
\begin{prf}
That $\widetilde{P}$ satisfies the Fay-shuffle equations is seen by a direct computation. On the other hand, from the definition of $\FSh_2(N)$, an element $P(X,Y) \in \FSh_2(N) \setminus \FSh^{pol}_2(N)$ necessarily has the form
\begin{equation} \label{eqn:decomp}
a\frac{X^{n-1}}{Y}+b\frac{Y^{n-1}}{X}+Q(X,Y),
\end{equation}
where $Q(X,Y) \in \Q[X,Y]$ is some specific polynomial, and $a,b \in \Q$ are not both equal to zero. Hence, we have 
\begin{equation}
\dim_{\Q}\FSh_2(N)^{pol}+1 \leq \dim_{\Q}\FSh_2(N) \leq \dim_{\Q}\FSh_2(N)^{pol}+2.
\end{equation}
But from the shuffle equation $P(X,Y)+P(Y,X)=0$, one sees that in \eqref{eqn:decomp}, it holds that $a=-b$, and thus $\dim_{\Q}\FSh_2(N)=\dim_{\Q}\FSh_2(N)^{pol}+1$.
\end{prf}
In order to prove Theorem \ref{thm:Fayshuf}, it therefore suffices to compute the dimension of $\FSh_2(N)^{pol}$, in the case where $N$ is odd.\footnote{The following argument that simplifies the original proof was communicated to the author by Francis Brown.} For this, denote by $W_N \subset \Q[X,Y]$ the subspace of homogeneous polynomials of degree $N$, which satisfy
\begin{equation} \label{eqn:w}
P(X,Y)+P(Y,X)=0 \quad P(X,Y)+P(Y,-X-Y)+P(-X-Y,X)=0.
\end{equation}
By comparing these equations with \eqref{eqn:fsh}, one sees immediately that 
\begin{equation} \label{eqn:wfsh}
W_{N-2}=\FSh_2(N)^{pol}, \quad \mbox{if $N$ is odd and $N \geq 3$}.
\end{equation}
Moreover, $W_N$ has the following representation-theoretic interpretation. A polynomial $P \in W_N$ can be viewed as a representation $\calp$ of the symmetric group $S_3$: the two-dimensional $\Q$-vector space $\calp \subset \Q[X,Y]$ spanned by $P(X,Y),P(Y,-X-Y)$ is acted upon by $S_3$ by
\begin{equation}
(12) \mapsto (P(X,Y) \mapsto P(-Y,-X)), \quad (123) \mapsto (P(X,Y) \mapsto P(Y,-X-Y)).
\end{equation}
Note that $(12)$ acts with trace $0$ and that $(123)$ acts with trace $-1$. Thus, it follows from the elementary representation theory of the symmetric group $S_3$ (cf.~e.g.~\cite{Ser}, Chap.~2) that $\calp$ is necessarily isomorphic to the unique irreducible, two-dimensional $S_3$-representation $\calw$. Also, since $P$ has degree $N$, and $\Q[X,Y]$ is naturally isomorphic to $\bigoplus_{N \geq 0}\Sym^N \calw$, the representation $\calp$ can be viewed as a sub-representation of $\Sym^N \calw$. Conversely, given a sub-representation $\calp$ of $\Sym^N \calw$ isomorphic to $\calw$, one verifies directly that there exists a unique (up to a non-zero scalar) element $P \in \calp$ such that $P \in W_N$.
\begin{prop} \label{prop:brown}
We have
\begin{equation}
\dim_{\Q} W_N=\left\lfloor \frac{N+2}{3} \right\rfloor.
\end{equation}
\end{prop}
\begin{prf}
Consider the decomposition
\begin{equation}
\Sym^N \calw \cong \calw_1 \oplus \dots \oplus \calw_n
\end{equation}
into irreducible sub-representations. By the preceding discussion, $\dim_{\Q} W_N$ is precisely equal to the number of summands in the above decomposition of $\Sym^N \calw$, which are isomorphic to $\calw$. Therefore, using character theory (cf.~\cite{Ser}, Theorem~2.4), one sees that
\begin{align}
\sum_{N \geq 0}\dim_{\Q} W_N t^N&=\frac{1}{6}\left(2\sum_{N \geq 0}(\Tr((1))_{\vert \Sym^N \calw})t^N-2\sum_{N \geq 0}(\Tr((123))_{\vert \Sym^N \calw})t^N\right) \notag\\
&=\frac{1}{6}\left(2\frac{1}{(1-t)^2}-2\frac{1}{1+t+t^2}\right) \notag\\
&=\frac{t}{(1-t)^2(1+t+t^2)} \notag\\
&=\sum_{N \geq 0}\left(\left \lfloor \frac {N+2}{3} \right \rfloor \right)t^{N+1},
\end{align}
as desired.
\end{prf}
Theorem \ref{thm:Fayshuf} now follows from combining Proposition \ref{prop:decomp} and Proposition \ref{prop:brown}, using in addition \eqref{eqn:wfsh}.
\bigskip

\noindent
4.6. {\bf Linear independence of elliptic double zeta values.}
We use the differential equation for A-elliptic multiple zeta values (Theorem \ref{thm:diff}) to prove a linear independence result for elliptic double zeta values.
\begin{thm} \label{thm:linind}
Let $N > 0$ be odd, and $k=\lfloor \frac N3 \rfloor$.
\begin{enumerate}
\item
The set
\begin{equation}
\calb_{N,2}=\left\{\frac{\partial}{\partial \tau} I^A(r,N-r;\tau) \; \vert \; 0 \leq r \leq k \right\}
\end{equation}
is linearly independent over $\Q$.
\item
We have 
\begin{equation}
\dim_{\Q}\cL_2\eZ^A_N \geq \left\lfloor \frac N3 \right\rfloor+1.
\end{equation}
\end{enumerate}
\end{thm}
\begin{prf}
First note that i) implies ii), as the A-elliptic double zeta values $I^A(r,N-r)$ for $r=0,\dots,k$ are linearly independent, since, by i), their derivatives are.

So assume there exists a relation
\begin{align} \label{eqn:diffrel}
\sum_{r=0}^k\lambda_r\frac{\partial}{\partial \tau}I^A(r,N-r;\tau)=0.
\end{align}
with $\lambda_r \in \Q$ for $r=0,\dots,k$. Then
\begin{align} \label{eqn:prematrix}
0=\lambda_0(NG_{0}(\tau)I^A(N)&-NG_{N+1}(\tau)I^A(0))+\sum_{r=1}^k\lambda_r\bigg(-(N-r)G_{N-r+1}(\tau)I^A(r) \notag\\
&+rG_{r+1}(\tau)I^A(N-r)-(-1)^rNG_{N+1}(\tau)I^A(0) \notag\\
&+\sum_{s=1}^{N+1}(N-s)\left(\binom{s-1}{r-1}-\binom{s-1}{N-r-1}\right)G_{N+1-s}(\tau)I^A(s)\bigg),
\end{align}
by Corollary \ref{cor:diff}. Denote by $C_N$ the $(k+1) \times (N+3)/2$ matrix whose entry $(C_N)_{r,s}$ is given by
\begin{equation}
(C_N)_{r,s}=\mbox{Coefficient of $\lambda_rI^A(2s)G_{N+1-2s}(\tau)$ in \eqref{eqn:prematrix}}.
\end{equation}
Since the Eisenstein series are linearly independent over $\C$, the existence of rational numbers $\lambda_0,\dots,\lambda_k$ solving \eqref{eqn:diffrel} is equivalent to the row vector $\Lambda_N:=(\lambda_0,\dots,\lambda_k) \in \Q^{k+1}$ solving
\begin{equation} \label{eqn:mat}
\Lambda_N \cdot C_N=(0,\dots,0).
\end{equation}
Thus if we can prove that the rank of $C_N$ is at least $k+1$, we are done, because then \eqref{eqn:mat} has only the trivial solution. Also, note that the first row of $C_N$ is equal to $(-N,0,\dots,0,N)$ by \eqref{eqn:prematrix}. Therefore, if we can prove that among the columns of $C_N$ indexed by $s=1,\dots,(N-1)/2$, there are $k$ linearly independent ones, $C_N$ will have rank at least $k+1$. For this, we can clearly assume that $N \geq 3$.

To this end, consider, for $N \geq 3$ and $N$ odd, the square submatrix $C'_N$ of $C_N$ consisting of the columns $s=1,2,\dots,k$. Looking at \eqref{eqn:prematrix} and comparing with Corollary \ref{cor:diff}, we see that its entries are given by
\begin{align} 
c'_{N;i,j}=(N-2j-2)\left(\binom{2j+1}{i}-\delta_{2j+1,i}\right), \quad 0 \leq i,j \leq k-1.
\end{align}
Since $N-2j-2 \neq 0$ for every $j$, as $N$ is odd, it is enough to prove that the scaled matrix $M_{k-1}=(m_{i,j})$ with
\begin{equation}
m_{i,j}=\left(\binom{2j+1}{i}-\delta_{2j+1,i}\right), \quad 0 \leq i,j \leq k-1
\end{equation}
is invertible for every $k \geq 1$, which is proved in the appendix (Proposition \ref{prop:app}).
\end{prf}
\bigskip

\noindent
4.7. {\bf End of the proof of the main result.}
We can now give a proof of Theorem \ref{mthm}.
\begin{prf}
From Theorem \ref{thm:Fayshuf}, we know that, if $N$ is even 
\begin{equation}
\dim_{\Q} \gr^{\cL}_2\eZ^A_N \leq \dim_{\Q}\FSh_2(N)=0,
\end{equation}
which proves the theorem in that case.

For odd $N$, we know by Theorem \ref{thm:Fayshuf}, using that $\cL_1\eZ^A_N=0$ if $N$ is odd,
\begin{equation}
\dim_{\Q}\gr^{\cL}_2\eZ^A_N \leq \left\lfloor \frac N3 \right\rfloor+1.
\end{equation}
On the other hand, by Theorem \ref{thm:linind}, the elliptic double zeta values
\begin{equation}
I^A(r,N-r;\tau) \quad 0 \leq r \leq k,
\end{equation}
where $k=\lfloor \frac N3 \rfloor$, are linearly independent over $\Q$. Therefore we have
\begin{equation}
\dim_{\Q}\gr^{\cL}_2\eZ^A_N = \left\lfloor \frac N3 \right\rfloor+1,
\end{equation}
if $N$ is odd.
Finally, the map
\begin{align}
\FSh_2(N)^* &\rightarrow \gr^{\cL}_2\eZ^A_N \notag\\
(X^*)^{k-1}(Y^*)^{l-1} &\mapsto I^A(k,l;\tau) \mod \cL_1\eZ^A_N
\end{align}
is surjective by construction, hence an isomorphism, since both sides have the same dimension. Since there are no non-trivial $\Q$-linear relations between elliptic double zeta values of different weights (cf.~Theorem~\ref{thm:weight}), it follows moreover that all $\Q$-linear relations in $\gr^{\cL}_2\eZ^A_N$ are a consequence of Fay and shuffle relations.
\end{prf}

\section{Outlook}
It is a natural problem to extend the computation of $D_{N,l}=\dim_{\Q}\gr^{\cL}_l\eZ^A_N$ carried out in this paper for $l=1,2$ to higher lengths. A partial result in this direction is that for odd $N$, we have
\[
D_{N,3}=\left\lfloor \frac{N+1}{6} \right\rfloor=\frac{N+1}{2}-\left\lfloor \frac N3 \right\rfloor-1.
\]
This follows from Theorem \ref{mthm}, together with the fact that every product $\pi^{2r}\calg_{2s}(\tau)$,
where $\calg_{2s}(\tau)$ denotes the indefinite Eisenstein integral of weight $2s$ (cf.~Section~2.4), arises as an A-elliptic multiple zeta value of length three and weight $N=2r+2s-1$. The latter can be proved again using the differential equation (Theorem \ref{thm:diff}).

At the moment, there is no analogous result for $D_{N,3}$ if $N$ is even or for higher length in general. The main obstacle in applying the methods of the length two case is simply that both the formula for the differential equation, as well as the combinatorics of the Fay-shuffle space become much more complicated in higher lengths. However, using the explicit $q$-expansions of A-elliptic multiple zeta values (cf.~\cite{BMS}, Section~2.3), one can obtain lower bounds $D^{low}_{N,l}$ for $D_{N,l}$. For example, for $l=3$ one has

\begin{center}
\begin{tabular}{|c||c|c|c|c|c|c|c|c|c|c|c|}
\hline
$N$ &0 &2 &4 &6 &8 &10 &12 &14 &16 &18 &20\\
\hline \hline
$D^{low}_{N,3}$ &0 &2 &3 &5 &8 &11 &14 &19 &23 &28 &34\\
\hline
\end{tabular}
\end{center}
It would be desirable to find a (conjectural) closed formula for $D_{N,3}$, and eventually, of course, for $D_{N,l}$ in general. Since elliptic multiple zeta values are related to the Lie algebra $\fu^{geom}$ of special derivations on the free Lie algebra on the set $\{a,b\}$ \cite{BMS}, it is possible that such a closed formula can be obtained from the dimensions of $\fu^{geom}$, which are explicitly known in depths $1,2$ and $3$, due to Brown \cite{BrownLetter}.
\begin{appendix}

\section{A binomial determinant} \label{ssec:bindet}
The proof of Theorem \ref{thm:linind}, and thus of the main theorem \ref{mthm} depended on the invertibility of a matrix with entries given by certain binomial coefficients. In this appendix, we complete the proof of Theorem \ref{thm:linind} by computing the determinant of the aforementioned binomial matrix explicitly.

Let $n$ be a positive integer and consider the matrix $M_n=(m_{i,j})_{0 \leq i,j \leq n}$ with
\begin{equation}
m_{i,j}=\binom{2j+1}{i}-\delta_{2j+1,i},
\end{equation}
where $\delta$ denotes the Kronecker delta. The proof of Theorem \ref{thm:linind} depends on the invertibility of $M_n$, which is proved in the next
\begin{prop} \label{prop:app}
We have
\begin{equation}
\det(M_n)=(2n+1)!!=1\cdot 3\cdot 5\dots\cdot(2n+1).
\end{equation}
In particular, $M_n$ is invertible for every $n$.
\end{prop}
The idea of the proof of Proposition \ref{prop:app} is to find a suitable LU-decomposition for $M_n$. We first need a lemma about binomial coefficients.
\begin{lem} \label{lem:app}
For all $a,b \geq 0$, we have
\begin{equation} \label{eqn:binom1}
\binom{a}{b}=\sum_{k=0}^b\binom{a-b+k}{k}\binom{a-b+1}{a-2b+2k+1}.
\end{equation}
\end{lem}
\begin{prf}
We first assume that $b\leq a/2$. In that case, the right hand side is equal to
\begin{equation} \label{eqn:app1}
\binom{a-b+1}{a-2b+1}{}_{3}F_{2}\left[\begin{smallmatrix}a-b+1,&-b/2+1/2,&-b/2\\\\&a/2-b+1,&a/2-b+3/2\end{smallmatrix};1\right],
\end{equation}
where ${}_{3}F_{2}$ is a hypergeometric function (cf.~\cite{B}, Chapter~II). If $b$ is even, we can apply Saalsch"utz's Theorem about ${}_3F_2$ (cf.~\cite{B}, Section~2.2), namely
\begin{equation}
{}_3F_2\left[\begin{smallmatrix}a,&b,&-n\\\\&c,&a+b+c+1-n\end{smallmatrix};1\right]=\frac{(c-a)_n(c-b)_n}{(c)_n(c-a-b)_n},
\end{equation}
where $(m)_k:=m(m+1)\dots(m+k-1)$ denotes the Pochhammer symbol. Hence, \eqref{eqn:app1} is equal to
\begin{equation}
\binom{a-b+1}{a-2b+1}\frac{(-a/2)_{b/2}(a/2-b/2+1/2)_{b/2}}{(a/2-b+1)_{b/2}(-a/2+b/2-1/2)_{b/2}}=\binom{a}{b}.
\end{equation}
In the case where $b$ is odd, we can apply Saalsch"utz's Theorem after interchanging $-b/2$ and $-b/2+1/2$ in the argument of ${}_3F_2$ in \eqref{eqn:app1} above, and get the result.
If $b>a/2$, then the same argument as above works, with $b$ replaced by $a-b$ throughout.
\end{prf}
Now consider the matrices $L_n=(l_{i,j})_{0 \leq i,j \leq n}$, $U_n=(u_{i,j})_{0 \leq i,j \leq n}$ with
\begin{equation}
l_{i,j}=\binom{j}{i-j}, \quad u_{i,j}=\begin{cases}1,&i=0\\ \binom{2j-i}{i-1}\frac{2j+1}{i},&0<i<2j+1\\0,& i \geq 2j+1\end{cases}.
\end{equation}
$L_n$ is a lower triangular matrix with determinant $1$, while $U_n$ is an upper triangular matrix with determinant $(2n+1)!!$. Hence, Proposition \ref{prop:app} follows from
\begin{lem}
For every $n$, we have
\begin{equation} \label{eqn:LU}
M_n=L_n U_n.
\end{equation}
\end{lem}
\begin{prf}
The assertion of the lemma is equivalent to
\begin{align} \label{eqn:binom2}
\binom{2j+1}{i}-\delta_{2j+1,i}=\sum_{k=1}^nl_{i,k}u_{k,j},
\end{align}
for all $i,j$ such that $0 \leq i,j \leq n$. For $i \geq 2j+1$, both sides of \eqref{eqn:binom2} are evidently equal to zero. If $i=0$, then, since $l_{0,k}=\delta_{0,k}$, both sides of \eqref{eqn:binom2} are equal to $1$.

It remains to prove \eqref{eqn:binom2} for $i,j$ such that $0<i<2j+1 \leq n$, i.e.
\begin{align} \label{eqn:binom3}
\binom{2j+1}{i}-\delta_{2j+1,i}&=\sum_{k=1}^n\binom{k}{i-k}\binom{2j-k}{k-1}\frac{2j+1}{k} \notag\\
&=\sum_{k=1}^i\binom{k}{i-k}\binom{2j-k}{k-1}\frac{2j+1}{k}.
\end{align}
Note that $\binom{k}{i-k}$ vanishes for $k > i$. We now rewrite the right hand side of \eqref{eqn:binom3} as
\begin{align}
\sum_{k=1}^i\binom{k}{i-k}\binom{2j-k}{k-1}\frac{2j+1}{k}&=\sum_{k=1}^i\binom{k}{i-k}\binom{2j+1-k}{k}\frac{2j+1}{2j+1-k} \notag \\
&=\sum_{k=0}^i\binom{2j+1-k}{i-k}\binom{2j+1-i}{2k-i}\frac{2j+1}{2j+1-k} \notag \\
&=\sum_{k=0}^i\binom{2j+1-k}{2j+1-i}\binom{2j+1-i}{2j-2k+1}\frac{2j+1}{2j+1-k} \notag \\
&=\sum_{k=0}^i\binom{2j-k}{2j-i}\binom{2j+1-i}{2j-2k+1}\frac{2j+1}{2j+1-i} \notag \\
&=\sum_{k=0}^i\binom{2j-k}{i-k}\binom{2j+1-i}{2j-2k+1}\frac{2j+1}{2j+1-i} \notag \\
&=\sum_{k=0}^i\binom{2j-i+k}{k}\binom{2j+1-i}{2j-2i+2k+1}\frac{2j+1}{2j+1-i},
\end{align}
and apply Lemma \ref{lem:app} to conclude the proof.
\end{prf}

\end{appendix}

\begin{bibtex}[\jobname]

@ARTICLE{BMS,
   author = {{Broedel}, J. and {Matthes}, N. and {Schlotterer}, O.},
    title = "{Relations between elliptic multiple zeta values and a special derivation algebra}",
  journal = {ArXiv e-prints},
  volume   = {hep-th/1507.02254}
archivePrefix = "arXiv",
   eprint = {1507.02254},
 primaryClass = "hep-th",
 keywords = {High Energy Physics - Theory, Mathematics - Number Theory},
     year = 2015,
    month = jul,
   adsurl = {http://adsabs.harvard.edu/abs/2015arXiv150702254B},
  adsnote = {Provided by the SAO/NASA Astrophysics Data System}
}

@ARTICLE{BMMS,
   author = {{Broedel}, J. and {Mafra}, C.~R. and {Matthes}, N. and {Schlotterer}, O.
	},
    title = "{Elliptic multiple zeta values and one-loop superstring amplitudes}",
  journal = {Journal of High Energy Physics},
archivePrefix = "arXiv",
   eprint = {1412.5535},
 primaryClass = "hep-th",
 keywords = {Scattering Amplitudes, Superstrings and Heterotic Strings, Conformal Field Models in String Theory},
     year = 2015,
    month = jul,
   volume = 7,
    pages = {112},
      doi = {10.1007/JHEP07(2015)112},
   adsurl = {http://adsabs.harvard.edu/abs/2015JHEP...07..112B},
  adsnote = {Provided by the SAO/NASA Astrophysics Data System}
}

@article{BrownLetter,
    author     =     {Brown, Francis},
    title     =     {{Letter to the author}},
    year     =     {2015},
    }

@MastersThesis{Pol,
    author     =     {Pollack, Aaron},
    title     =     {{Relations between derivations arising from modular forms}},
    school     =     {Duke University},
    year     =     {2009},
    }

@incollection {CEE,
             AUTHOR = {Calaque, Damien and Enriquez, Benjamin and Etingof, Pavel},
              TITLE = {Universal {KZB} equations: the elliptic case},
          BOOKTITLE = {Algebra, arithmetic, and geometry: in honor of {Y}u. {I}.
                       {M}anin. {V}ol. {I}},
             SERIES = {Progr. Math.},
             VOLUME = {269},
              PAGES = {165--266},
          PUBLISHER = {Birkh\"auser Boston, Inc., Boston, MA},
               YEAR = {2009},
            MRCLASS = {32G34 (11F55 17B37 20C08 32C38)},
           MRNUMBER = {2641173 (2011k:32018)},
         MRREVIEWER = {Gwyn Bellamy},
                DOI = {10.1007/978-0-8176-4745-2_5},
                URL = {http://dx.doi.org/10.1007/978-0-8176-4745-2\_5},
}
 
@article {Levin,
    AUTHOR = {Levin, Andrey},
     TITLE = {Elliptic polylogarithms: an analytic theory},
   JOURNAL = {Compositio Math.},
  FJOURNAL = {Compositio Mathematica},
    VOLUME = {106},
      YEAR = {1997},
    NUMBER = {3},
     PAGES = {267--282},
      ISSN = {0010-437X},
     CODEN = {CMPMAF},
   MRCLASS = {11F37 (11F27 11G40 11R70 19F27)},
  MRNUMBER = {1457106 (98d:11048)},
MRREVIEWER = {Alexey A. Panchishkin},
       DOI = {10.1023/A:1000193320513},
       URL = {http://dx.doi.org/10.1023/A:1000193320513},
}

@article{BL,
     author         = {Brown, Francis. and Levin, Andrey},
     title          = {Multiple elliptic polylogarithms},
	 journal 		= {ArXiv e-prints},
	 volume   		= {math.NT/1110.6917},
     year           = {2010},
     eprint         = {1110.6917},
     archivePrefix  = {arXiv},
     primaryClass   = {math},
}

@book {W,
    AUTHOR = {Weil, Andr{\'e}},
     TITLE = {Elliptic functions according to {E}isenstein and {K}ronecker},
      NOTE = {Ergebnisse der Mathematik und ihrer Grenzgebiete, Band 88},
 PUBLISHER = {Springer-Verlag, Berlin-New York},
      YEAR = {1976},
     PAGES = {ii+93},
      ISBN = {3-540-07422-8},
   MRCLASS = {10DXX (01A55 10-03)},
  MRNUMBER = {0562289 (58 \#27769a)},
MRREVIEWER = {S. Chowla},
}

@unpublished{BrownAnatomy,
     author         = "Brown, F.",
     title          = "{Anatomy of an associator}",
     year			= "{2012}",
     note			= "{Notes}"
}

@article{LR,
     author         = {A. Levin and G. Racinet},
     title          = {Towards multiple elliptic polylogarithms},
	 journal 		= {ArXiv e-prints},
	 volume   		= {math.NT/0703237v1},
     year           = {2007},
     eprint         = {0703237v1},
     archivePrefix  = {arXiv},
     primaryClass   = {math},
}

@article{E,
      hyphenation   = {american},
      author        = {Benjamin Enriquez},
      title         = {Analogues elliptiques des nombres multiz{\'e}tas},
      journal 		= {ArXiv e-prints},
	  volume   		= {math.NT/1301.3042},
      date          = {2013-01-14},
      year          = {2013},
      eprinttype    = {arxiv},
      archivePrefix = {arXiv},
      eprint        = {1301.3042}
} 

@article{MMV,
     author         = {Brown, Francis},
     title          = {Multiple modular values for ${\SL_2(\Z)}$},
     year           = {2014},
     journal		= {ArXiv e-prints},
     volume			= {math.NT/1407.5167v1},
}

@article {E2,
    AUTHOR = {Enriquez, Benjamin},
     TITLE = {Elliptic associators},
   JOURNAL = {Selecta Math. (N.S.)},
  FJOURNAL = {Selecta Mathematica. New Series},
    VOLUME = {20},
      YEAR = {2014},
    NUMBER = {2},
     PAGES = {491--584},
      ISSN = {1022-1824},
   MRCLASS = {17B35 (11M32 14H10 16S30 20F36)},
  MRNUMBER = {3177926},
       DOI = {10.1007/s00029-013-0137-3},
       URL = {http://dx.doi.org/10.1007/s00029-013-0137-3},
}

@article {Z,
    AUTHOR = {Zagier, Don},
     TITLE = {Periods of modular forms and {J}acobi theta functions},
   JOURNAL = {Invent. Math.},
  FJOURNAL = {Inventiones Mathematicae},
    VOLUME = {104},
      YEAR = {1991},
    NUMBER = {3},
     PAGES = {449--465},
      ISSN = {0020-9910},
     CODEN = {INVMBH},
   MRCLASS = {11F67 (11F27 11F55)},
  MRNUMBER = {1106744 (92e:11052)},
MRREVIEWER = {Rolf Berndt},
       DOI = {10.1007/BF01245085},
       URL = {http://dx.doi.org/10.1007/BF01245085},
}

@article {LM,
    AUTHOR = {Le, Thang Tu Quoc and Murakami, Jun},
     TITLE = {Kontsevich's integral for the {K}auffman polynomial},
   JOURNAL = {Nagoya Math. J.},
  FJOURNAL = {Nagoya Mathematical Journal},
    VOLUME = {142},
      YEAR = {1996},
     PAGES = {39--65},
      ISSN = {0027-7630},
     CODEN = {NGMJA2},
   MRCLASS = {57M25 (11M99)},
  MRNUMBER = {1399467 (97d:57009)},
MRREVIEWER = {Sergei K. Lando},
       URL = {http://projecteuclid.org/euclid.nmj/1118772043},
}

@preamble{
   "\def\cprime{$'$} "
}
@article {Dr,
    AUTHOR = {Drinfel{\cprime}d, V. G.},
     TITLE = {On quasitriangular quasi-{H}opf algebras and on a group that
              is closely connected with {${\rm Gal}(\overline{\bf Q}/{\bf
              Q})$}},
   JOURNAL = {Algebra i Analiz},
  FJOURNAL = {Algebra i Analiz},
    VOLUME = {2},
      YEAR = {1990},
    NUMBER = {4},
     PAGES = {149--181},
      ISSN = {0234-0852},
   MRCLASS = {16W30 (17B37)},
  MRNUMBER = {1080203 (92f:16047)},
MRREVIEWER = {Ivan Penkov},
}

@incollection {Man,
    AUTHOR = {Manin, Yuri I.},
     TITLE = {Iterated integrals of modular forms and noncommutative modular
              symbols},
 BOOKTITLE = {Algebraic geometry and number theory},
    SERIES = {Progr. Math.},
    VOLUME = {253},
     PAGES = {565--597},
 PUBLISHER = {Birkh\"auser Boston, Boston, MA},
      YEAR = {2006},
   MRCLASS = {11F67 (11G55 11M41)},
  MRNUMBER = {2263200 (2008a:11062)},
MRREVIEWER = {Caterina Consani},
       DOI = {10.1007/978-0-8176-4532-8_10},
       URL = {http://dx.doi.org/10.1007/978-0-8176-4532-8\_10},
}

@article {BK,
    AUTHOR = {Broadhurst, D. J. and Kreimer, D.},
     TITLE = {Association of multiple zeta values with positive knots via
              {F}eynman diagrams up to {$9$} loops},
   JOURNAL = {Phys. Lett. B},
  FJOURNAL = {Physics Letters. B},
    VOLUME = {393},
      YEAR = {1997},
    NUMBER = {3-4},
     PAGES = {403--412},
      ISSN = {0370-2693},
     CODEN = {PYLBAJ},
   MRCLASS = {11M41 (11Z05 57M25 81T18)},
  MRNUMBER = {1435933 (98g:11101)},
MRREVIEWER = {Louis H. Kauffman},
       DOI = {10.1016/S0370-2693(96)01623-1},
       URL = {http://dx.doi.org/10.1016/S0370-2693(96)01623-1},
}

@incollection {Del,
   AUTHOR = {Deligne, P.},
    TITLE = {Le groupe fondamental de la droite projective moins trois
             points},
BOOKTITLE = {Galois groups over ${\bf Q}$ ({B}erkeley, {CA}, 1987)},
   SERIES = {Math. Sci. Res. Inst. Publ.},
   VOLUME = {16},
    PAGES = {79--297},
PUBLISHER = {Springer, New York},
     YEAR = {1989},
  MRCLASS = {14G25 (11G35 11M06 11R70 14F35 19E99 19F27)},
 MRNUMBER = {1012168 (90m:14016)},
MRREVIEWER = {James Milne},
      DOI = {10.1007/978-1-4613-9649-9_3},
      URL = {http://dx.doi.org/10.1007/978-1-4613-9649-9\_3},
}

@article {IKZ,
    AUTHOR = {Ihara, Kentaro and Kaneko, Masanobu and Zagier, Don},
     TITLE = {Derivation and double shuffle relations for multiple zeta
              values},
   JOURNAL = {Compos. Math.},
  FJOURNAL = {Compositio Mathematica},
    VOLUME = {142},
      YEAR = {2006},
    NUMBER = {2},
     PAGES = {307--338},
      ISSN = {0010-437X},
   MRCLASS = {11M41},
  MRNUMBER = {2218898 (2007e:11110)},
MRREVIEWER = {David Bradley},
       DOI = {10.1112/S0010437X0500182X},
       URL = {http://dx.doi.org/10.1112/S0010437X0500182X},
}

@book {Bou,
    AUTHOR = {Bourbaki, N.},
     TITLE = {\'{E}l\'ements de math\'ematique. {F}asc. {XXXIV}. {G}roupes
              et alg\`ebres de {L}ie. {C}hapitre {IV}: {G}roupes de
              {C}oxeter et syst\`emes de {T}its. {C}hapitre {V}: {G}roupes
              engendr\'es par des r\'eflexions. {C}hapitre {VI}: syst\`emes
              de racines},
    SERIES = {Actualit\'es Scientifiques et Industrielles, No. 1337},
 PUBLISHER = {Hermann, Paris},
      YEAR = {1968},
     PAGES = {288 pp. (loose errata)},
   MRCLASS = {22.50 (17.00)},
  MRNUMBER = {0240238 (39 \#1590)},
MRREVIEWER = {G. B. Seligman},
}

@book {B,
    AUTHOR = {Bailey, W. N.},
     TITLE = {Generalized hypergeometric series},
    SERIES = {Cambridge Tracts in Mathematics and Mathematical Physics, No.
              32},
 PUBLISHER = {Stechert-Hafner, Inc., New York},
      YEAR = {1964},
     PAGES = {v+108},
   MRCLASS = {33.20 (40.00)},
  MRNUMBER = {0185155 (32 \#2625)},
}

@article {Ch,
    AUTHOR = {Chen, Kuo Tsai},
     TITLE = {Iterated path integrals},
   JOURNAL = {Bull. Amer. Math. Soc.},
  FJOURNAL = {Bulletin of the American Mathematical Society},
    VOLUME = {83},
      YEAR = {1977},
    NUMBER = {5},
     PAGES = {831--879},
      ISSN = {0002-9904},
   MRCLASS = {55D35 (58A99)},
  MRNUMBER = {0454968 (56 \#13210)},
MRREVIEWER = {Jean-Michel Lemaire},
}

@book {Ser,
    AUTHOR = {Serre, Jean-Pierre},
     TITLE = {Repr\'esentations lin\'eaires des groupes finis},
   EDITION = {revised},
 PUBLISHER = {Hermann, Paris},
      YEAR = {1978},
     PAGES = {182},
      ISBN = {2-7056-5630-8},
   MRCLASS = {20-01 (20C99)},
  MRNUMBER = {543841 (80f:20001)},
}
\end{bibtex}
\bibliographystyle{abbrv}
\bibliography{\jobname}

\bigskip

\noindent
{\small
{\it E-mail: }\texttt{nils.matthes@uni-hamburg.de}
\bigskip

\noindent {\sc Fachbereich Mathematik (AZ)\\ Universit\"at Hamburg\\ Bundesstrasse 55\\ D-20146 Hamburg}}

\end{document}

%% file: template.tex
\usepackage[T1]{fontenc}
\usepackage[english]{babel}
\usepackage[latin1]{inputenc}
\usepackage{lmodern}
\usepackage{ngerman}
\usepackage{amsmath}
\usepackage{amssymb}
\usepackage{amscd}
\usepackage{enumitem}
\usepackage[bookmarksopen]{hyperref}
\usepackage{amsthm}
\usepackage{dsfont}
\usepackage{makeidx}
\usepackage{fancyhdr}
\usepackage{mathrsfs}
\usepackage{tikz}
\usepackage{pstricks}
\usepackage{graphicx}
\usepackage{color}
\usepackage{ifpdf}
\usepackage{wrapfig}
\usepackage{currfile}
\usepackage{sectsty}
\usepackage[OT2,T1]{fontenc}
\usepackage{array}
\usepackage{shuffle}
\usepackage{titlesec}
\DeclareGraphicsExtensions{.eps,.bmp}

\usetikzlibrary{backgrounds}

\theoremstyle{definition}
\newtheorem{dfn}{Definition}[section]
\newtheorem{rmk}[dfn]{Remark}

\theoremstyle{plain}
\newtheorem{prop}[dfn]{Proposition}

\newtheorem{conj}[dfn]{Conjecture}
\newtheorem{thm}[dfn]{Theorem}
\newtheorem{lem}[dfn]{Lemma}
\newtheorem{cor}[dfn]{Corollary}
\newtheorem*{mthm}{Main Theorem}

\newenvironment{prf}{\begin{proof}[{\it Proof: \nopunct}]}{\end{proof}}

\def\H          {\mathbb H}
\def\N          {\mathbb N}
\def\R          {\mathbb R}
\def\Q          {\mathbb Q}

\def\C          {\mathbb C}
\def\Z          {\mathbb Z}

\def\calo        {\mathcal O}

\def\cE            {\mathcal E}
\def\eZ           {\mathcal{EZ}}
\def\Eis           {\mathcal{E}\mathrm{is}}
\def\calb         {{\mathcal B}}
\def\cD           {{\mathcal D}}

\def\calg         {{\mathcal G}}
\def\cali         { {\mathcal I}}

\def\cL         {\mathcal{L}}

\def\calp         {\mathcal{P}}

\def\calw         {\mathcal W}

\def\cZ			  {\mathcal{Z}}

\def\fu           {\mathfrak{u}}

\def\be                 {\begin{equation}}
\def\ee                 {\end{equation}}

\def\Im         {\mathrm{Im}}

\def\gr       {\mathrm{gr}}
\def\dd       {\mathrm{d}}

\def\Reg      {\mathrm{Reg}}

\DeclareSymbolFont{cyrletters}{OT2}{wncyr}{m}{n}
\DeclareMathSymbol{\sha}{\mathalpha}{cyrletters}{"58}

\DeclareMathOperator{\SL}{SL}

\DeclareMathOperator{\ad}{ad}

\DeclareMathOperator{\FSh}{FSh}
\DeclareMathOperator{\Sym}{Sym}
\DeclareMathOperator{\Tr}{Tr}

\setenumerate{label=\textnormal{(\roman*)}, ref=(\textnormal{\roman*})}

\sectionfont{\normalfont\centering\fontsize{15}{18}\selectfont}
\subsectionfont{\fontsize{11}{11}\selectfont}

\titleformat{\subsection}[runin]
  {\normalfont\bfseries}{\thesubsection}{0.5em}{}

\renewcommand{\thesubsection}{\thesection.\arabic{subsection}.}

\numberwithin{equation}{section}